\documentclass[11pt]{article}
\usepackage[cp866]{inputenc}
\begin{document}

\centerline{\large\bf Internal connection between the field theory equations.}
\centerline{\large\bf Fundamentals of the field theory} 
\centerline {\bf L.I. Petrova}
\centerline{\it Moscow State University, Department  of Computational Mathematics and}
\centerline{\it Cybernetics, Russia, e-mail: ptr@cs.msu.su}

\renewcommand{\abstractname}{Abstract}
\begin{abstract}

It is shown that there is a correspondence between field theory equations such as the Dirac, Shr\H{o}dinger, Maxwell, Einstein equations and closed exterior forms of a certain degree. 

In this case, the Dirac and Shr\H{o}dinger equations for the wave function correspond to closed exterior forms of zero degree. The Shr\H{o}dinger equation for the state functional corresponds to closed exterior forms of the first degree. The Maxwell's equations based on exterior forms of second degree. Einstein's equation for the gravitational field consists of covariant tensors, which correspond to closed exterior forms of the second degree. However, the covariant tensors of the Einstein equation are derived from the covariant tensors, which correspond to closed exterior forms of third degree.

Such a correspondence between the field theory equations and closed exterior forms of a certain degree reveals the internal connection between the field theory equations. 

At the same time, it was shown that closed exterior forms, on which the field theory equations for physical fields are based, are associated with the equations of mathematical physics for material media, such as thermodynamic, gas-dynamic, electromagnetic, cosmological equations, etc. 

This indicates the connection between the field theory equations and the mathematical physics equations and reveals the foundations of field theory.

\end{abstract}
 
 
\section{Introduction}
 
As is known, field theory equations describing physical fields, such as the equations of Dirac, Shr\H{o}dinger, Maxwell, Einstein, are based on the invariance and covariance properties. 

It turns out that there is a mathematical formalism that possesses these properties. This is the theory of closed skew-symmetric differential forms [1,2,3].  

Closed exterior forms are defined on integrable manifolds and satisfy the condition: $d\omega \,=0$. They are differentials, and therefore invariants, and correspond to the conservation laws for physical fields: the presence of conservative quantities or structures. (The example is Noether's theorem.) Moreover, the dual forms corresponding to closed inexact (defined only on integrable structures) exterior forms have covariant properties. 

It is shown that there is a correspondence between the field theory equations and closed exterior forms. This reveals an internal connection between theories describing different types of physical fields. 

In addition, it is shown that closed exterior forms, on which the field theory equations for physical fields are based, are connected with the mathematical physics equations for material media, such as thermodynamic, gas-dynamic, electromagnetic, cosmological equations, etc. 

This reveals the basis of field theory equations and the nature of physical fields. 

\section{Correspondence between field theory equations and closed exterior forms of a certain degree}

The study of field theory equations by closed external skew-symmetric forms showed that there is a correspondence between the field theory equations and closed external skew-symmetric forms. It can be seen that the field theory equations are associated with closed exterior forms of a certain degree. 

Closed exterior forms of zero degree correspond to quantum mechanics equations.

The Dirac bra- and cket- vectors form a closed exterior form of zero degree.

The Shr\H{o}dinger equation for the wave function:
$$
\Huge i\hbar\frac{\partial\Psi}{\partial t}=\hat{H}\Psi
$$ 
consists of operators corresponding to closed exterior forms of zero degree.

The Heisenberg equation for quantum mechanics is associated with a dual form corresponding to a closed exterior form of zero degree.
 
\bigskip
The Hamiltonian formalism is based on the properties of a closed exterior form of the first degree and the corresponding dual form. The Shr\H{o}dinger equation for the action functional is such an equation: 
$$
h{{\partial \psi }\over {\partial t}}\,=\,H\psi
$$
Here the Hamilton function $H$ is a closed skew-symmetric form of the first degree.  

\bigskip
The properties of closed exterior forms of the second degree and their dual forms underlie the electromagnetic field equations. Covariant electromagnetic field tensor $F_{\mu\nu}$ in  Maxwell's equations obeys the identical relations $d\theta^2=0$, $d^*\theta^2=0$, where $\theta^2=\frac{1}{2}F_{\mu\nu}dx^\mu dx^\nu$ is a closed exterior form of the second degree, and $^*\theta^2$ is the corresponding dual form.
                    
\bigskip

Einstein's equation 
$$G_{\mu\nu}={{8{\pi}G}\over{c^4}}T_{\mu\nu}$$
(where $G_{\mu\nu}=R_{\mu\nu}-1/2g_{\mu\nu}R$ is the Einstein tensor) is a relation in the covariant tensors, which are related to closed exterior and dual forms of the second degree. In this case, the Einstein equation is obtained from the relation with covariant tensors, which are analogues of closed exterior forms of the third degree. 

\subsection*{Approach to the unified field theory.} 

It can be seen that the field theory equations are associated with closed exterior forms of a certain degree. This shows that there is an internal connection between field theory equations that describe different types of physical fields.

Obviously, that the degree of closed exterior forms is a parameter that unites field theories into the unified field theory. 

In addition, it can be shown that the degree of closed exterior forms is a parameter that unifies interactions. If we denote the degree of closed exterior forms by $k$, then $k=0$ corresponds to the strong interaction, $k=1$ corresponds to weak interaction, $k=2$ corresponds to electromagnetic, and $k=3$ corresponds to the gravitational interaction [3, Appendix 5]. 

One more general property of field theory equations can be noted here. They are relations for functionals such as the wave function, the action functional, the Einstein tensor. For the electromagnetic field such a functional is the Poynting vector [4]. 

Thus, one can see that the field theory equations are relations for functionals and are based on the properties of closed exterior forms.

\section{What is the physical meaning of closed exterior forms and how are they being obtained?}

It turns out that closed exterior forms, on which the field theory equations that describe physical fields are based, are connected with the equations of mathematical physics that describe material media, such as thermodynamic, gas-dynamic, electromagnetic, cosmological, etc. 

And this is due to the conservation laws.

The conservation laws have a feature.

As is known, conservation laws for physical fields are the presence of conservative quantities or structures. Closed exterior forms describe such conservative quantities or structures (as already noted, Noether's theorem is the example). 

And the conservation laws for material media are the conservation laws of energy, momentum, angular momentum and mass. These are differential conservation laws. They are described by differential equations. Equations of mathematical physics describing material media, as is known, consist of such equations. 

From the equations of conservation laws for material media, skew-symmetric forms also are obtained. But these are skew-symmetric forms that are evolutionary, since they are defined, in contrast to exterior forms, on a nonintegrable manifold. 

It turns out that such evolutionary skew-symmetric forms, which are obtained from the mathematical physics equations, generate closed exterior forms that underlie the field theory equations. 

Such results were obtained in the study of the consistency of the conservation laws equations for material media. (It should be noted that the quantum properties of the mathematical physics equations and the field theory equations expressly depend on the consistency of the conservation laws equations.) 

\subsection{Investigation of the conservation laws equations consistency for material media.}

Some features of the study should be noted. 

Conservation law equations were expressed in terms of state functionals [4], such as wave function, entropy, action functional, etc., which are also functionals of the field theory equations. 

In addition, the conservation laws equations were written in frame of reference, associated with the studied material medium. This is an accompanying frame of reference is connected with the manifold constructed of 
the trajectories of the material medium elements. 
(Lagrange coordinate system is  example of such frame of reference.) 

\bigskip
Let us investigate the energy and momentum conservation laws equations consistency. 

In accompanying frame of reference,  the energy conservation law equation is written as:
$$
{{\partial \psi }\over {\partial \xi ^1}}\,=\,A_1 \eqno(1)
$$
Here $\xi^1$ - coordinate along the trajectory, $A_1$ is a quantity that depends on the material medium properties and on external energy actions onto the medium.
 
Similarly, in accompanying frame of reference, the momentum conservation law equation reduces to the equation:
$$
{{\partial \psi}\over {\partial \xi^{\nu }}}\,=\,A_{\nu },\quad \nu \,=\,2,\,...\eqno(2)
$$ 
where $\xi ^{\nu }$ are coordinates in the normal direction to the trajectory, $A_{\nu }$ - quantities that depend on the medium properties and external force actions.

Equations (1) and (2) can be reduced to the relation:
$$
d\psi\,=\,\omega \eqno(3)
$$
where
$\omega \,=\,A_{\mu }\,d\xi ^{\mu }$ is a skew-symmetric differential form of the first degree. (The summation is carried out over the repeating index.) 

In the general case (for the conservation equations of energy, momentum, angular momentum and mass), this relation will have the form:
$$
d\psi \,=\,\omega^p \eqno(4)
$$
where $\omega^p$ is a skew-symmetric differential form of the degree $p=0,1,2,3$.  

Since the conservation laws equations are evolutionary, the skew-symmetric form and relations (3) and (4) are also evolutionary.

\subsection{Features of the evolutionary skew-symmetric form and evolutionary relation.}

The evolutionary form, in contrast to the exterior form, is defined on an deformable, nonintegrable manifold (constructed of the trajectories of the material medium elements). 
The differential of a form whose basis is a nonintegrable manifold is not equal to zero. 

This can be shown on the example of the evolutionary form of the first degree 
$\omega \,=\,A_{\alpha }\,d\xi ^{\alpha }$.   

Differential $d\omega$  of the evolutionary form 
$\omega \,=\,A_{\alpha }\,d\xi ^{\alpha }$ can be written as 
$d\omega=K_{\alpha\beta}d\xi^\alpha d\xi^\beta$, where  
$K_{\alpha\beta}=A_{\beta;\alpha}-A_{\alpha;\beta}$ is a commutator of the evolutionary form 
$\omega$, and  
$A_{\beta;\alpha}$ , $A_{\alpha;\beta}$  are covariant derivatives. If we express the covariant derivatives in terms of connectedness, then the form $\omega$ commutator can be written as: 
$$
K_{\alpha\beta}=\left(\frac{\partial A_\beta}{\partial
\xi^\alpha}-\frac{\partial A_\alpha}{\partial
\xi^\beta}\right)+(\Gamma^\sigma_{\beta\alpha}-
\Gamma^\sigma_{\alpha\beta})A_\sigma \eqno(5)
$$
The evolutionary form commutator is not equal to zero, since it includes connectedness of the metric form of a nonintegrable manifold that are not symmetric, i.e. the expression 
$(\Gamma^\sigma_{\beta\alpha}-\Gamma^\sigma_{\alpha\beta})$ 
not equal to zero [5]. 
                               
Since the commutator is not equal to zero, then, consequently, the differential of the evolutionary form is also not equal to zero. And this means that the evolutionary form, unlike the exterior form, cannot be a closed form, i.e. a differential. 

It can be shown that the evolutionary form  $\omega^p$ is not a differential. 

(The example of such an evolutionary form is the Lagrangian.)

\bigskip 
Since on the right side of the evolutionary relation $d\psi \,=\,\omega^p$  there is an evolutionary form that is not a differential, this means that the evolutionary relation turns out to be a nonidentical relation (there is a differential on the left, and an evolutionary form that is not a differential on the right).

(An example of a nonidentical evolutionary relation is the first law of thermodynamics.) 

{\footnotesize [It can be noted that the nonidentity of the evolutionary relation, which was obtained in the study of the conservation laws equations consistency, indicates that on the original coordinate space the conservation laws equations for material media turn out to be inconsistent, i.e., the conservation laws are noncommutative. The hidden properties of the mathematical physics equations [6] are connected precisely with the noncommutativity of the conservation laws.]}   

\section{Correspondence between field theory equations and evolutionary relation}

As noted above, the field theory equations are relations for the $\psi$ functionals. 

The resulting evolutionary relation $d\psi \,=\,\omega^p$ is a relation for the same functionals.
  
However, these are different relations.

The field theory equations are identical relations, since the right side of the relations contains closed exterior forms, i.e. differentials, as in the left side. 

And the evolutionary relation is a nonidentical relation, since the evolutionary form on the right side of the evolutionary relation is not a closed exterior form, i.e., a differential.

But it turns out that identical relations, which are analogous to the field theory equations, are obtained from the evolutionary relation.

\subsection{Realization of the identity relation.}

An identity relation can be realized from a nonidentical evolutionary relation if a closed exterior form, which is a differential, is obtained from an evolutionary form that is not closed.

But this can only happen under a degenerate transformation, i.e., under a transformation that does not preserve the differential, since the differential of the evolutionary form is not equal to zero, while the differential of a closed exterior form is equal to zero. 

\subsection*{The unique meaning of the degenerate transformation.} 

In field theory, as is known, nondegenerate transformations are used, i.e., transformations that preserve the differential. They make transitions between integrable structures.

A degenerate transformation is a transformation that does not preserve the differential. It has a unique property: under a degenerate transformation, an integrable structure arises and a transition occurs from a nonintegrable space to structures of an integrable manifold. 

An example of a degenerate transformation is the Legendre transformation. It is the Legendre transformation that makes the transition from a Lagrangian, nonintegrable manifold to an  Hamiltonian structures of integrable manifold.

Degenerate transformations occur discretely, arise spontaneously when any degrees of freedom are realized. This corresponds to turning to zero of such functional expressions as Jacobians, determinants, Poisson brackets, residues, etc. Such conditions describe a realized integrable structures.

\subsection*{Degenerate transformation. Realization of a closed external inexact skew-symmetric form.} 

A degenerate transformation is a transformation of skew-symmetric differential forms. It is described by the transition from the evolutionary form differential, which is not equal to zero, to the differentials of the closed inexact exterior form and the corresponding dual form, which are equal to zero.

This can be written as:

$d\omega^p\ne 0 \to $ (degenerate transformation) $\to d_\pi{}^*\omega^p=0$, $ d_\pi \omega^p=0$

The $d_\pi{}^*\omega^p=0$ condition indicates the implementation of a closed dual form $^*\omega^p$, which describes some realized integrable structure $\pi$ (pseudostructure by its metric properties). 

And the $ d_\pi \omega^p=0$ condition indicates that the closed inexact exterior form $\omega^p_\pi$ is realized on the integrable structure $\pi$.

\subsection*{Identical relation.}

Realization of a closed inexact exterior form $\omega^p_\pi$ and dual form $^*\omega^p$ leads to the fact that on the pseudostructure $\pi$ from the nonidentical evolutionary relation $d\psi \,=\,\omega^p$ the relation 
$$
d_\pi{\psi}=\omega^p_\pi\eqno(6) 
$$ 
which is the identical relation, is realized, since the right side has a closed form, which is a differential, like the left side. 

\section{Identical relations of the mathematical physics equations as an analogue of the field theory equations}

As has been shown, the field theory equations describing physical fields are relations for field theory functionals and are based on the closed exterior forms properties. 

Identical relations  $d_\pi{\psi}=\omega^p_\pi$ derived from the mathematical physics equations, describing material media, are relations for the same functionals and contain closed exterior forms. 

It turns out that there is a connection between the closed exterior forms of the field theory equations and the closed exterior forms of the identical relations obtained from the mathematical physics equations. The right parts of the field theory equations are obtained from the mathematical physics equations under a degenerate transformation.

Thus, the Hamilton function in the Hamiltonian formalism is obtained from the Euler equation for the Lagrange function under the Legendre transformation, which is a degenerate transformation [7]. 

The covariant tensors in the Maxwell's equations are connected to the Poynting vector, which describes the energy flow of electromagnetic media. 

The energy-momentum tensor in the Einstein equation follows from the analysis of the equations of the conservation laws of energy, momentum and the equation of state [8], i.e. also from the mathematical physics equations. 
{\footnotesize [Here it should be noted that when deriving the energy-momentum tensor, Einstein imposed a number of assumptions on the conservation laws equations in order to prove the covariance of the energy-momentum tensor [8]. However the covariance of the energy-momentum tensor is satisfied under a degenerate transformation, in which the consistency of the energy and momentum equations is realized.]} 

The connection of the field theory equations with the identical relations of the mathematical physics equations reveals the foundations of the field theory and can serve as an approach to the general field theory. 

\subsection*{Quantum nature of the field theory equations.} 

The correspondence between the field theory equations and the identical relations of the mathematical physics equations reveals the quantum character of the field theory equations. 

As has been shown, the identical relations of mathematical physics are realized discretely under a degenerate transformation.

The connection of the field theory equations with the identical relations of mathematical physics equations indicates that the field theory equations also realize discretely. In this case, a structure appears (which is described by a dual form) with a conservative quantity (a closed exterior form), i.e., an invariant object, a quantum, arises.

This reveals the quantum character of the field theory equations.

In the quantum mechanics, this is taken into account by the discrete Hamiltonian spectrum [7].

In the quantum physics, this is connected to the discrete properties of the Hamilton function, which is defined on sections of the bundle of the cotangent space [7].

The correspondence between the field theory equations and the identical relations of the mathematical physics equations also substantiates the Einstein equation discreteness.
When deriving the Einstein equation [8], it was assumed that the transitions from the curvature tensor to Einstein's tensor (at which the Bianchi identity is fulfilled, the connectedness  coefficients are symmetric ones, i.e., the connectedness coefficients are the Christoffel symbols, and there exists a transformation under which the connectedness coefficient becomes zero [9]) are nondegenerate transformations (transitions were made from the tensors upper indices to the lower indices, which is true only for tensors defined on integrable manifolds).
However, on spaces with curvature, which are non-integrable manifolds, such transformations are nondegenerate, and are realized only discretely, when any degrees of freedom are realized. 

It turns out that the Einstein equation can only be satisfied discretely [9]. In this case, it should be taken into account that quantum effects arise in discrete transitions. (In the material media, they appear as observable formations, such as waves, vortices, pulsations, etc.)

\bigskip
Here you can pay attention to the following. The covariance conditions, which the field theory equations must satisfy, are satisfied only for symmetric connections. Symmetric connections are obtained by a degenerate transformation from nonsymmetric connections.

Einstein, as shown in the book [10], in the last years of his life tried to obtain symmetric connections from asymmetric ones in order to derive the general field theory equation. He, as noted by Pauli [11], tried to find a nondegenerate transformation that would make it possible to transform nonsymmetric connections into symmetric ones. But such attempts did not lead to success, because the transformation that transforms nonsymmetric connections into symmetric ones is a degenerate transformation. 

\section{Conclution}

It has been shown that there is a mathematical formalism that reveals the properties of the field theory equations. These are skew-symmetric differential forms. 

It is shown that field theory equations are based on the properties of closed exterior forms. At the same time, field theory equations describing physical fields of various types are associated with closed exterior forms of a certain, corresponding degree. 

This reveals the internal connection of field theory equations, unifies field theories, and can serve as an approach to a unified field theory. 

\bigskip

It turns out that from the mathematical physics equations describing material media, evolutionary skew-symmetric forms are obtained, which generate closed exterior forms that underlie the field theory equations.

This points to the connection of the field theory equations with the mathematical physics equations and reveals the foundations and nature of the field theory equations.

[1] {\it Cartan E.,} Les Systems Differentials Exterieus ef Leurs Application
Geometriques. -Paris, Hermann, 1945.

[2] {\it Petrova L.I.,} Exterior and evolutionary differential forms in mathematical physics: Theory and Applications, -Lulu.com,157, 2008. 

[3] Petrova L.I., Role of skew-symmetric differential forms in mathematics, 2010.  
https://arxiv.org/pdf/1007.4757.pdf 

[4]  L.I. Petrova, Physical meaning and a duality of concepts of wave function, action functional, entropy, the Pointing vector, the Einstein tensor. https://arxiv.org/pdf/1001.1710.pdf 

[5] {\it Tonnelat M.-A.,} Les principles de la theorie electromagnetique et la relativite, Masson, Paris, 1959.

[6]  Petrova Ludmila, Discrete Quantum Transitions, Duality: Emergence of Physical Structures and Occurrence of Observed Formations (Hidden Properties of Mathematical Physics Equations). //Applied Mathematics and Physics,  Scientific Research Publishing (United States), vol 8, No. 9, 1911-1921, 2020. 

[7] Petrova L.,  Qualitative Investigation of Hamiltonian Systems by Application of Skew-Symmetric Differential Forms, //Symmetry, MDPI (Basel, Switzerland), Vol 13, No.1, 1-7, 2021. 

[8] {\it A. Einstein,} The Meaning of Relativity. Princeton, 1953.

[9]  Petrova L., Evolutionary Relation of Mathematical Physics Equations. Evolutionary Relation as Foundation of Field Theory. Interpretation of the Einstein Equation. //Axioms, MDPI (Basel, Switzerland), Vol 10(2), No. 46, 1-10, 2021. 
https://doi.org/10.3390/axioms10020046 

[10] {\it Abraham PAIS,} THE SCIENCE AND THE LIFE OF Albert EINSTEIN, Oxford univerity New York Toronto Melbourno,1982.
 
[11] {\it W. Pauli,} Theory of Relativity, Pergamon Press, 1958. 

\end{document}